\documentclass[12pt]{amsart}

\usepackage[backref=page]{hyperref}

\usepackage{amssymb, amsthm,hyperref,xcolor, mathrsfs, cancel, fullpage, comment, thmtools}

\usepackage[capitalise]{cleveref}

\newtheorem{theorem}{Theorem}[section]
\newtheorem{lemma}[theorem]{Lemma}

\theoremstyle{remark}

\numberwithin{equation}{section}

\newcommand{\mo}{{-1}}

\newcommand{\aux}{\ensuremath{\mathrm{aux}}}

\title{On the cop number and the weak Meyniel conjecture for algebraic graphs}

\author{Arindam Biswas}

\address{Polynom, 20 Rue Jacques Daguerre, Rueil-Malmaison, Paris, 92500}
\email{arin.math@gmail.com}
\thanks{}

\author{Jyoti Prakash Saha}
\address{Department of Mathematics, Indian Institute of Science Education and Research Bhopal, Bhopal Bypass Road, Bhauri, Bhopal 462066, Madhya Pradesh, India}
\curraddr{}
\email{jpsaha@iiserb.ac.in}
\thanks{}

\subjclass[2010]{05C57, 91A24, 91A43}

\keywords{Cops and robbers, Weak Meyniel conjecture, Cayley sum graphs, generalised Cayley graphs}

\begin{document}

\maketitle

\begin{abstract}
We show that the cop number of the Cayley sum graph of a finite group $G$ with respect to a symmetric subset $S$ is at most twice its degree when the graph is connected, undirected. We also prove that a similar bound holds for the cop number of generalised Cayley graphs and twisted Cayley sum graphs under some conditions. These extend a result of Frankl to such graphs. Using the above bounds and a result of Bollob\'{a}s--Janson--Riordan, we show that the weak Meyniel conjecture holds for these algebraic graphs.
\end{abstract}

\section{Introduction}

The game of cops and robbers was introduced by Nowakowski and Winkler \cite{NowakowskiWinklerVertexToVertexPursuitInGraph}, and independently by Quilliot \cite{QuilliotThesis}, though other variants of the game, similar in style have been played before mainly using Alquerque boards. As mentioned by Bonato, a reference to the game of cops and robbers exists in a book of puzzles by Dudeney 
\cite{DudeneyAmusement}, as early as 1917.
This game is played on a finite connected graph as follows. First, the cops each choose a vertex in the given graph. Then, the robber chooses a vertex. The game begins and the cops and the robber alternate in taking turns by moving along the edges, with the cops starting first. 
The minimum number of cops, required to capture the robber no matter how the robber plays, is called the \emph{cop number} of the graph. For further details, we refer to the book by Bonato--Nowakowski \cite{BonatoNowakowskiSTMLCopsRobbers}, and the survey by Baird--Bonato \cite{BairdBonatoMeynielConjSurvey}.

Since its inception, the game has been studied by numerous authors including Aigner--Fromme \cite{AignerFrommeGameOfCopsRobbers}, Andreae \cite{AndreaeNoteOnPursuitGamePlayedGraphs, AndreaePursuitGameGraphsMinorExcluded}, Hamidoune \cite{HamidounePursuitGameOnCayleyGraphs}, Maamoun--Meyniel \cite{MaamounMeynielGameOfPolicemenRobber}, Nowakowski--Winkler \cite{NowakowskiWinklerVertexToVertexPursuitInGraph}, Quilliot \cite{QuilliotThesis,QuilliotDiscretePursuitGame} etc. Frankl, in \cite{FranklPursuitGameOnCayleyGraph}, showed that for the Cayley graph of a finite, abelian group with respect to a generating set of size $d$, its cop number is at most $\lceil \frac{d+1}{2} \rceil$. One of the most important conjecture in this area is Meyniel's conjecture. Henri Meyniel, in a personal communication with Frankl, conjectured that the cop number of a connected graph on $n$ vertices is $O(\sqrt n)$ \cite{FranklCopsRobbersGraphsLargeGirthCayleyGraphs}. 
Upper bounds on the cop number were obtained by Frankl \cite{FranklCopsRobbersGraphsLargeGirthCayleyGraphs}, Chiniforooshan \cite{ChiniforooshanBetterBddForCopNumberGeneralGraph}, with further improvements by Lu--Peng \cite{LuPengMeynielConjOfTheCopNumber}, 
Frieze--Krivelevich--Loh
\cite{FriezeKrivelevichLohVariationsOnCopsRobbers}, 
Scott--Sudakov \cite{ScottSudakovBddForTheCopsRobbersProblem}. Bollob\'{a}s--Kun--Leader \cite{BollobasKunLeaderCopsRobbersInRandomGraph} 
showed that the cop number of the binomial random graph $G(n,p)$ is $O(\sqrt n \log n)$ asymptotically almost surely if $p \geq (2 + \varepsilon) \log n /n$ for some $\varepsilon > 0$. 
Soon after, the $\log n$ factor was removed by Pra{\l}at--Wormald \cite{PralatWormaldMeynielConjRandomGraphs} for $p\geq (\frac 12 + \varepsilon) \log n /n$ for some $\varepsilon > 0$. They also obtained the same bound for sparse random graphs \cite{PralatWormaldMeynielConjRandomDReguGraphs}. Recently, Bradshaw \cite{BradshawProofOfMeynielConjAbelianCayleyGraphs} and Bradshaw--Hosseini--Turcotte \cite{BradshawHosseiniTurcotteCopsRobbersDirUnditAbelianCayley} proved Meyniel's conjecture for abelian Cayley graphs. 
For further details, we refer to the works 
\cite{LuPengMeynielConjOfTheCopNumber, LehnerCopNumberOfToroidalGraphs, HosseiniMoharGonzalezMeynielConjGraphsBddDegre, GonzalezSebastianHosseiniKnoxMoharReedCopsRobberOriToroiGrid, DasGahlawatSahooSenCopsRobberSomeFamiliesOriGra, BradshawHosseiniMoharStachoOnTheCopNumberGraphWithHighGirth, PralatWhenDoesARandomGraphHaveConsCopNumber, BradshawHosseiniTurcotteCopsRobbersDirUnditAbelianCayley, HasiriShinkarMeynielExtremalFamiliesAbelianCayley} among others. 
For general classes of graphs, this conjecture is still wide open. Indeed, even the weak Meyniel conjecture is open, which states that for a fixed constant $c>0$, the cop number of a graph on $n$ vertices is $O(n^{1-c})$. In this article, one of our goals is to show the weak Meyniel conjecture for various classes of algebraic graphs.

\subsection{Algebraic graphs} 
Throughout the article, $G$ denotes a finite group, and $S$ denotes a subset of $G$. The Cayley graph $C(G,S)$ is the graph having the elements of $G$ as vertices, and there is an edge from a vertex $u$ to a vertex $v$ if $v=u s$ for some $s\in S$, while the Cayley sum graph $C_{\Sigma}(G,S)$   has the elements of $G$ as vertices, and there is an  edge from a vertex $u$ to a vertex $v$ if $v=u^{-1} s$ for some $s\in S$. There are important structural differences between Cayley graphs and Cayley sum graphs. 
For example, the former is vertex-transitive, while the latter may fail to be so. 

Cayley sum graphs are interesting objects of study in their own right which is apparent from a flurry of works on the subject including those of Chung \cite{Chung89JAMS}, and more recently those of Green \cite{GreenCountingSetsWithSmallSubsetsClique, GreenChromatic}, Grynkiewicz--Lev--Serra \cite{GrynkiewiczLevSerraConnCaylSum}, Green--Morris \cite{GreenMorrisCountingSetsWithSmallSubsets}, DeVos--Goddyn--Mohar--\v{S}\'{a}mal \cite{DeVosGoddynMoharSamalCayleySumFullerene}, Lev \cite{LevSumDiffHamiltCycle}, Amooshahi--Taeri \cite{AmooshahiTaeri}, the authors \cite{CheegerCayleySum, VertexTra} etc. However, even though the game of cops and robbers on Cayley graphs has an extensive literature since the time of Frankl, the game has never been studied on Cayley sum graphs before (to the best of the knowledge of the authors). In this article, we address this issue. 
We also study it on twisted Cayley graphs and twisted Cayley sum graphs. 
For a discussion on these graphs, we refer to \cite{CheegerTwisted, Expansion}. 
The notion of twisted Cayley graphs is related to the notion of generalised Cayley graphs as defined by Maru\v{s}i\v{c}, Scapellato and Zagaglia Salvi \cite{MarusicScapellatoZagagliaSalviGeneralizedCayleyGraph}. 

\subsection{Results obtained}
Under certain assumptions, we show that for any of the above-mentioned graphs, the cop number is at most twice their degree. This extends some previous results of Frankl \cite[Theorem 1.5]{FranklCopsRobbersGraphsLargeGirthCayleyGraphs}. 
See Theorems 
\ref{Thm:CayleySumCopNoBdd}, \ref{Thm:CayleyTwistedCopNoBdd}, \ref{Thm:CayleyTwistedSumCopNoBdd} 
for the precise statements. Combining these bounds with a result of Bollob\'{a}s--Janson--Riordan, we show under suitable hypotheses that a stronger version of the weak Meyniel conjecture (in fact, Meyniel's conjecture up to a logarithmic factor) holds for these graphs of algebraic nature (\cref{Thm:WMCIntro}), to be referred to as \emph{algebraic graphs} for brevity.

We assume throughout that the robber makes a move along some edge during each round. If this hypothesis is relaxed, then employing an additional cop, who keeps the distance with the robber equal or lower when the robber makes a move, and reduces the distance with the robber at least by $1$ when the robber stays idle, it follows that the cop number would increase at most by $1$. 
We remark that the remaining cops will move as per the strategy (to be 
described during the proofs) whenever the robber makes a move, 
and would remain idle if the robber does not move. 

\begin{theorem}[Cop number of undirected Cayley sum graphs]
\label{Thm:CayleySumCopNoBdd}
Suppose the Cayley sum graph $C_\Sigma(G, S)$ is undirected and connected. 
Assume that $S$ is symmetric. Then the cop number of $C_\Sigma(G, S)$ is at most $2|S|$. 
\end{theorem}

For an automorphism $\sigma$ of $G$, the twisted Cayley graph $C(G, S)^\sigma$ has $G$ as its set of vertices, and there is an edge from $x$ to $y$ if $y = \sigma(xs)$ for some $s\in S$.

\begin{theorem}[Cop number of undirected twisted Cayley  graphs]
\label{Thm:CayleyTwistedCopNoBdd}
Let $\sigma$ be an automorphism of $G$ of order two. 
Suppose the twisted Cayley graph $C(G, S)^\sigma$ is undirected and connected. 
Assume that $S$ is closed under conjugation by the elements of $G$ and $S$ is symmetric. Then the cop number of $C(G, S)^\sigma$ is at most $2|S|$. 
\end{theorem}

For an automorphism $\sigma$ of $G$, the twisted Cayley sum graph $C_\Sigma(G, S)^\sigma$ has $G$ as its set of vertices, and there is an edge from $x$ to $y$ if $y = \sigma(x^\mo s)$ for some $s\in S$.

\begin{theorem}[Cop number of undirected twisted Cayley sum graphs]
\label{Thm:CayleyTwistedSumCopNoBdd}
Let $\sigma$ be an automorphism of $G$ of order two. 
Suppose the twisted Cayley sum graph $C_\Sigma(G, S)^\sigma$ is undirected and connected. 
Assume that $S$ is symmetric. Then the cop number of $C_\Sigma (G, S)^\sigma$ is at most $2|S|$. 
\end{theorem}

\begin{theorem}
[\cref{Thm:WMC}]
\label{Thm:WMCIntro}
The cop number of full Cayley graphs, and the algebraic graphs considered in Theorems 
\ref{Thm:CayleySumCopNoBdd}, \ref{Thm:CayleyTwistedCopNoBdd}, \ref{Thm:CayleyTwistedSumCopNoBdd}
is at most $2\sqrt n \log n$, where $n$ denotes the cardinality of the underlying group. 
In particular, the weak Meyniel conjecture holds for these graphs. 
\end{theorem}

\section{Strategy for bounding the cop number}

In the following, we adopt the convention  that during the 
zeroth round of the cops and robber game on a graph, the cops and the robber occupy their initial positions. 
For $n\geq 1$, during 
the $n$-th round of the cops and robber game on a graph, the robber takes its $n$-th turn and the cops take their $n$-th turn. 
Let $\Gamma$ be a graph among the graphs $C_\Sigma(G, S), C(G, S)^\sigma, C_\Sigma(G, S)^\sigma$ as considered in Theorems 
\ref{Thm:CayleySumCopNoBdd}, \ref{Thm:CayleyTwistedCopNoBdd}, \ref{Thm:CayleyTwistedSumCopNoBdd}. 
From the hypothesis of any of these results, it follows that $S$ generates $G$. 

We provide an outline of the general framework to bound the cop number of $C_\Sigma(G, S)$, $C(G, S)^\sigma$, $C_\Sigma(G, S)^\sigma$. 
The precise details of the strategy depend on these graphs, which are to be specified in the subsequent sections.

\subsection{The moves of the primary cops}

We introduce $|S|$ \emph{primary} cops and let them occupy the vertex at the identity element of $G$, and let the robber occupy some vertex of $\Gamma$. Henceforth, the primary cops are to be referred to as cops, for brevity.

\subsubsection{Notations}

At the end of the zeroth round, that is, after the cops and the robber have occupied their initial positions on $\Gamma$, we label the cops using the elements of $S$, and the cop corresponding to an element $s$ of $S$ is denoted by $c(s)$. 
The cops will be relabelled at the end of each round. 
Note that  a cop may have different labels during various rounds. 
For instance, a cop at the beginning of a round may have the label $c(s_1)$, at the beginning of the next round, the same cop may have the label $c(s_2)$ for $s_1, s_2\in S$. 
We describe the move of the cops after a move of the robber, and the method of relabelling of the cops.

Let $n$ be a positive integer. 
Suppose at the beginning of the $n$-th round, the robber is at the vertex $y$ and moves to $\tilde y$, and a cop $C$ having label $c(s)$ is positioned at the vertex $z$. The \emph{connection} of the robber, is denoted by $\xi_{y, \widetilde y}$, and to be defined as an element of $S$ depending on $\Gamma$ and the elements $y,\tilde y$. 
The \emph{gap} between the cop $C$ and the robber at the beginning of the $n$-th round, is denoted by $\Delta_{z, y}$, and to be defined as an element of $G$ depending on $\Gamma$. Since $S$ generates $G$, we may write the gap as $\Delta_{z, y}= w_{z, y,k}s^i$ where $i$ is an integer, $w_{z, y, k}$ denotes the product of (a choice of) $k$ elements of $S$ with $k$ as small as possible. The nonnegative integer $k$ is defined as the \emph{tail} of the cop $C$ with label $c(s)$ at the beginning of the $n$-th round. If the cop $C$ with label $c(s)$ has tail zero at the beginning of the $n$-th round, then the \emph{power} of $C$ is defined as the smallest non-negative integer $i$ such that $\Delta_{z, y}= s^i$.

\subsubsection{The strategy}
To describe the move of the cop $C$ with label $c(s)$, and the method of relabelling $C$, we will consider the situations when 
\begin{enumerate}
\item 
the connection of the robber is not a power of the label $s$ of $C$,
\item 
the connection of the robber is a power of the label $s$ of $C$,
\begin{enumerate}
\item 
the tail $k$ of the cop $C$ is positive with $w_{z, y, k} = g_1 g_2 \cdots g_k$ for some $g_1, g_2, \ldots , g_k\in S$, 

\item 
the tail of $C$ is zero,

\begin{enumerate}
\item the connection of the robber does not lie in $\{s, s^\mo\}$, 

\item the connection of the robber lies in $\{s, s^\mo\}$.
\end{enumerate}
\end{enumerate}
\end{enumerate}
The above will be referred to as Case (1), (2), (2a), (2b), (2bi) etc. 
In the subsequent sections, depending on the graph $\Gamma$, the strategy for $C$ is to be described considering the above cases. 

\subsubsection{Reduction of tail and power}

In the subsequent sections, for each of the graphs $C_\Sigma(G, S)$, $C(G, S)^\sigma$, $C_\Sigma(G, S)^\sigma$, after describing the strategy for the cops, we will establish the following lemma.

\begin{lemma}
\label{Lemma:ReductionGen}
Let $n$ be a positive integer, and $C$ be a primary cop. 
\begin{enumerate}

\item 
The tail of $C$ at the beginning of the $n$-th round is greater than or equal to its tail at the beginning of the $(n+1)$-st round. 
\item 
If the connection of the robber in the $n$-th round is a power of the label that the cop $C$ has at the beginning of the $n$-th round and 
the tail of the cop $C$ at the beginning of this round is $k$ for some integer $k\geq 1$, then the tail of the cop $C$ at the beginning of the $(n+1)$-st round is at most $k-1$. 
\item

If $C$ has tail zero and power $i$ at the beginning of the $n$-th round, then the power of $C$ at the beginning of the $(n+1)$-st round is 
\begin{enumerate}
\item 
$i$ if the connection of the robber during the $n$-th round is not the inverse of the label of $C$, 
\item 
$i-2$ if the connection of the robber during the $n$-th round is the inverse of the label of $C$ and $i \geq 2$.
\end{enumerate}
\end{enumerate}
\end{lemma}

Let us explain briefly the role of \cref{Lemma:ReductionGen} in capturing the robber. 
During any round, the connection of the robber is an element of $S$ and hence, it is the label of some cop at the beginning of that round. If the tail of that cop is positive in that round, then by \cref{Lemma:ReductionGen}(2), its tail strictly decreases. Further, during any round, the tail of none of the cops increases by \cref{Lemma:ReductionGen}(1). 
Hence, by the pigeonhole principle, after enough rounds, the tail of some cop $C$ reduces to zero, and it remains so in the next rounds too. 
Next, during enough of the upcoming rounds, if the connection of the robber coincides with the inverse of the label of some cop $D$ with tail zero, then the power of $D$ reduces by two in those rounds by \cref{Lemma:ReductionGen}(3b), and hence the power of $D$ would eventually reduce to zero (in which case, the robber would be captured by $D$), or to one (in which case, the capturing would be done by one of the additional \emph{auxiliary} cops, to be introduced later). 
By applying the pigeonhole principle more carefully, we would argue that one can find such a cop $D$.

\section{Cop numbers of Cayley sum graphs}
In this section, we bound the cop number of the graph $\Gamma = C_\Sigma(G, S)$ as in \cref{Thm:CayleySumCopNoBdd}. 
Note that the graph $C_\Sigma(G, S)$ is undirected if and only if $S$ is closed under conjugation by the elements of $G$. In the proof of \cref{Thm:CayleySumCopNoBdd}, the assumption that $S$ is symmetric is used crucially, for instance, to make sure that the proposed strategy for the cops consists of moves that are possible in $C_\Sigma(G, S)$. Moreover, if $S$ is symmetric, then it follows that the graph $C_\Sigma(G, S)$ is connected if and only if $S$ generates $G$.

Note that \cref{Thm:CayleyTwistedSumCopNoBdd} is more general, and \cref{Thm:CayleySumCopNoBdd} is a simpler (untwisted) version of \cref{Thm:CayleyTwistedSumCopNoBdd}. First, we establish \cref{Thm:CayleySumCopNoBdd} in order to explain the strategy of the proof and the argument in this easier case. Next, \cref{Thm:CayleyTwistedSumCopNoBdd} is established using an analogous strategy. 

\subsection{The moves of the primary cops}

Let $n$ be a positive integer. Suppose at the beginning of the $n$-th round, the robber is at the vertex $y$ and moves to $\tilde y = y^\mo t$ for some $t\in S$, and a cop $C$ having label $c(s)$ is positioned at the vertex $z$. The \emph{connection} of the robber is defined as 
$$
\xi_{y, \widetilde y} 
: = 
\begin{cases}
y^\mo t^\mo y & \text{ if $n$ is odd}, \\
t & \text{ if  $n$ is even}. 
\end{cases}
$$
The \emph{gap} between the cop $C$ and the robber at the beginning of the $n$-th round is defined as 
$$
\Delta_{z, y}:= 
\begin{cases}
z^\mo y & \text{ if $n$ is odd},\\
zy^\mo & \text{ if $n$ is even}.
\end{cases}
$$
We describe the move of the cop $C$ with label $c(s)$, and the method of relabelling $C$. 
If $n$ is odd, then 
$C$ moves to 
$$
\begin{cases}
\xi_{y, \widetilde y}^\mo z^\mo & \text{ when Case (1) holds}, \\
g_1^\mo z^\mo & \text{ when Case (2a) holds}, \\
\xi_{y, \widetilde y}^\mo z^\mo & \text{ when Case (2bi) holds}, \\
s^\mo z^\mo & \text{ when Case (2bii) holds}.
\end{cases}
$$
Note that the above moves are possible in $C_\Sigma(G, S)$ since for any $x\in S$, there is an edge from $z$ to any vertex of the form $x^\mo z^\mo$ because $x^\mo z^\mo = z^\mo (z x^\mo z^\mo)$ and $S$ is closed under conjugation and $S$ is symmetric. If $n$ is even, then $C$ moves to 
$$
\begin{cases}
z^\mo t& \text{ when Case (1) holds}, \\
z^\mo g_1& \text{ when Case (2a) holds}, \\
z^\mo t& \text{ when Case (2bi) holds}, \\
z^\mo s& \text{ when Case (2bii) holds}.
\end{cases}
$$
The label of $C$ is changed to $c(\xi_{y, \widetilde y}^\mo s \xi_{y, \widetilde y})$. Since $S$ is closed under conjugation, the element $\xi_{y, \widetilde y}^\mo s \xi_{y, \widetilde y}$ lies in $S$ and hence, the relabelling is justified. 

\subsubsection{Reduction of tail and power}
We prove \cref{Lemma:ReductionGen} for $\Gamma = C_\Sigma(G, S)$, showing that under the above strategy, the tail and the power of the cops do not increase, and in certain situations, they strictly decrease.

\begin{proof}
[Proof of \cref{Lemma:ReductionGen} for $\Gamma = C_\Sigma(G, S)$]
Assume that at the beginning of the $n$-th round, the robber is at the vertex $y$, and a cop $C$ is at the vertex $z$ with label $c(s)$, and the robber moves to $\tilde y = y^\mo t$. Suppose the cop $C$ moves to $\tilde z$, and its label changes to $c(\tilde s)$.

Suppose Case (1) holds, that is, the connection of the robber in the $n$-th round is not a power of the label that the cop $C$ has at the beginning of the $n$-th round. 
If $n$ is odd, then writing $z^\mo y = w_{z, y, k} s^i$, we note that the gap between the cop $C$ and the robber at the beginning of the $(n+1)$-st round satisfies 
\begin{align*}
\Delta_{\widetilde z, \widetilde y}
& = 
\tilde z \tilde y^\mo
\\
& = (\xi_{y, \widetilde y}^\mo z^\mo) (y^\mo t)^\mo \\
& = \xi_{y, \widetilde y}^\mo z^\mo y (y^\mo t^\mo y)\\
& = \xi_{y, \widetilde y}^\mo 
w_{z, y, k} s^i
\xi_{y, \widetilde y} \\
& = (\xi_{y, \widetilde y}^\mo 
w_{z, y, k}\xi_{y, \widetilde y})
(\xi_{y, \widetilde y}^\mo  s
\xi_{y, \widetilde y} )^i \\
& = (\xi_{y, \widetilde y}^\mo 
w_{z, y, k}\xi_{y, \widetilde y})
\tilde s^i .
\end{align*}
Since $S$ is closed under conjugation and $w_{z, y, k}$ 
is the product of some $k$ elements of $S$, it follows that 
$\xi_{y, \widetilde y}^\mo 
w_{z, y, k}\xi_{y, \widetilde y}$ is 
the product of some $k$ elements of $S$. 
If $n$ is even, then writing $z y^\mo = w_{z, y, k} s^i$, we note that the gap between the cop $C$ and the robber at the beginning of the $(n+1)$-st round satisfies 
\begin{align*}
\Delta_{\widetilde z, \widetilde y}
& = 
\tilde z^\mo \tilde y
\\
& = (z^\mo t )^\mo (y^\mo t) \\
& =  t^\mo zy^\mo t \\
& =  t^\mo  w_{z, y, k} s^i t \\
& =  t^\mo  w_{z, y, k} t (t^\mo  s t )^i \\
& =  (t^\mo  w_{z, y, k} t) \tilde s^i 
.
\end{align*}
Since $S$ is closed under conjugation, it follows that 
$t^\mo  w_{z, y, k} t$ is the product of some $k$ elements of $S$.
This proves the first statement of the lemma in Case (1). 

Suppose Case (2) holds, that is, the connection of the robber in the $n$-th round is a power of the label that the cop $C$ has at the beginning of the $n$-th round, that is, $\xi_{y, \widetilde y} = s^j$ for some integer $j$. 
Note that the label of $C$ remains unchanged in this round. 
First, let us consider Case (2a). So, the tail of the cop $C$ at the beginning of the $n$-th round is $k$ for some integer $k\geq 1$. 
If $n$ is odd, then writing $z^\mo y = g_1 g_2 \cdots g_k s^i$, we note that the gap between the cop $C$ and the robber at the beginning of the $(n+1)$-st round satisfies 
\begin{align*}
\Delta_{\widetilde z, \widetilde y}
& = 
\tilde z \tilde y^\mo
\\
& = g_1^\mo z^\mo  (y^\mo t)^\mo \\
& = g_2 \cdots g_k s^i y^\mo t^\mo y \\
& = g_2 \cdots g_k s^i \xi_{y, \widetilde y}\\
& = g_2 \cdots g_k s^{i+ j}
. 
\end{align*}
If $n$ is even, then writing $zy^\mo = g_1 g_2 \cdots g_k s^i$, we note that the gap between the cop $C$ and the robber at the beginning of the $(n+1)$-st round satisfies 
\begin{align*}
\Delta_{\widetilde z, \widetilde y}
& = 
\tilde z^\mo \tilde y
\\
& = (z^\mo g_1)^\mo y^\mo t\\
& = g_1^\mo z y^\mo t\\
& = g_2 \cdots g_k s^i  t\\
& = g_2 \cdots g_k s^i  \xi_{y, \widetilde y}\\
& = g_2 \cdots g_k s^{i + j}
.
\end{align*}
This proves the first statement in Case (2a), and also proves the second statement. 
Next, let us consider Case (2b). So, the tail of the cop $C$ at the beginning of the $n$-th round is zero. 
If $n$ is odd, then writing $z^\mo y = s^i$, we note that
\begin{align*}
\Delta_{\widetilde z, \widetilde y}
& = 
\tilde z \tilde y^\mo
\\
& = \xi_{y, \widetilde y}^\mo z^\mo   (y^\mo t)^\mo \\
& = \xi_{y, \widetilde y}^\mo z^\mo y  (y^\mo t^\mo y ) \\
& = \xi_{y, \widetilde y}^\mo s^i \xi_{y, \widetilde y}\\
& = s^i,
\end{align*}
and also note that 
\begin{align*}
\Delta_{\widetilde z, \widetilde y}
& = 
\tilde z \tilde y^\mo
\\
& = s^\mo z^\mo  (y^\mo t)^\mo \\
& = s^\mo s^i y^\mo  t^\mo y\\
& = s^{i-1} \xi_{y, \widetilde y},
\end{align*}
which is equal to $s^i$ (resp. $s^{i-2}$) when the connection $\xi_{y, \widetilde y}$ is equal to $s$ (resp. $s^\mo$). 
If $n$ is even, then writing $zy^\mo = s^i$, we note that 
\begin{align*}
\Delta_{\widetilde z, \widetilde y}
& = 
\tilde z^\mo \tilde y
\\
& = (z^\mo t)^\mo y^\mo t\\
& = t^\mo zy^\mo t\\
& = t^\mo s^{i} t\\
& = s^i, 
\end{align*}
and also note that 
\begin{align*}
\Delta_{\widetilde z, \widetilde y}
& = 
\tilde z^\mo \tilde y
\\
& = (z^\mo s)^\mo y^\mo t\\
& = s^\mo zy^\mo t\\
& = s^{i-1} t\\
& = s^{i-1} \xi_{y, \widetilde y},
\end{align*}
which is equal to $s^i$ (resp. $s^{i-2}$) when the connection $\xi_{y, \widetilde y}$ is equal to $s$ (resp. $s^\mo$). 
This proves the first statement in Case (2b), and also proves the third statement. This completes the proof of the lemma. 
\end{proof}

\subsection{The moves of the auxiliary cops} 
Let us introduce another set of $|S|$ \emph{auxiliary} cops and describe their moves. 
For each primary cop $C$ already employed, we introduce a cop $C_\aux$. During any round, $C_\aux$ occupies the vertex $v^\mo$ if $C$ is at the vertex $v$. Note that if $C$ moves from $u$ to $v$ during a round, then $uv$ lies in the symmetric set $S$, and hence there is an edge joining the vertex $v^\mo$ with $u^\mo$ in $C_\Sigma(G, S)$, and consequently, the cop $C_\aux$ can move from $u^\mo$ to $v^\mo$ in the undirected graph $C_\Sigma(G, S)$.

\subsection{Capturing the robber}
Using \cref{Lemma:ReductionGen}, we show that the cop number of $C_\Sigma(G, S)$ is at most $2|S|$. 

\begin{proof}
[Proof of \cref{Thm:CayleySumCopNoBdd}]
Let $x$ denote the position of the robber at the end of  the zeroth round. 
Recall that the cops occupy the vertex at the identity element of $G$ and so the gap between the cops and the robber at the beginning of the first round is $x$. 
Suppose $x$ can be expressed as a product of $m$ elements of $S$. 
By the pigeonhole principle, there is a cop $C$ such that within at most $(m + |G|)|S|$ rounds (these are successive and start from the first round), there are precisely $m+|G|$ rounds (henceforth called the \emph{special rounds}), such that for any of these $m+|G|$ rounds, the connection of the robber is equal to 
the label that the cop $C$ has at the beginning of that round. 
Note that if the labels of two cops at the beginning of a round are inverses of each other, then their labels at the beginning of the next round are also inverses of each other. In the following, $D$ denotes the cop whose label at the beginning of the first round is equal to the inverse of the label of the cop $C$. We remark that $C, D$ need not be different. Note that during any of the special rounds, the connection of the robber is equal to the inverse of the label that the cop $D$ has at the beginning. 
Since the cops occupy the vertex at the identity element of $G$ and $x$ can be expressed as a product of $m$ elements of $G$, it follows that the tail of every cop at the beginning of the first round is at most $m$. 
By \cref{Lemma:ReductionGen}, the tail of $C$ is zero during the $(m+1)$-st special round, and the next rounds. Similarly, the tail of $D$ also has this property. 
 
Suppose the power of  $D$ is equal  to $2\ell$ for some integer $0 \leq \ell \leq |G|/2$ at the beginning of the $(m+1)$-st special round. By \cref{Lemma:ReductionGen}, the power of   $D$ is equal  to zero at the end of the $(m+\ell)$-th special round, i.e., the  robber is captured by the cop $D$ at the end of the $(m+\ell)$-th special round. 
Suppose the power of  $D$ is equal  to $2\ell+1$ for some integer $\ell$  satisfying $1 \leq 2\ell+1 \leq |G|$ at the beginning of the $(m+1)$-st special round. By \cref{Lemma:ReductionGen}, the power of $D$ is equal to $1$ at the beginning of the $(m+\ell+1)$-st special round. Denote the position of the robber (resp. the cop $D$) at the beginning of the $(m+\ell+1)$-st special round by $y$ (resp. $z$). 
Let $s$ denote the label of $C$ at the beginning of this round. Suppose the robber moves to $y^\mo t$ during this round. 
If the $(m+\ell+1)$-st special round is an even round, then $zy^\mo = s^\mo$ and the connection of the robber is equal to $t$, which is equal to $s$, and consequently, $y^\mo t = y^\mo s = z^\mo$. 
If the $(m+\ell+1)$-st special round is an odd round, then $z^\mo y = s^\mo$ and the connection of the robber is equal to $y^\mo t^\mo y$, which is equal to $s$, and consequently, $y^\mo t = s^\mo y^\mo  = z^\mo$. 
Hence, in this special round, the robber is captured by the auxiliary cop $D_\aux$. Consequently, the robber will be captured within at most $(m + |G|/2)|S|$ rounds. 
\end{proof}

\section{Cop numbers of twisted Cayley graphs}
Let $\sigma$ and $S$ be as in \cref{Thm:CayleyTwistedCopNoBdd}. 
In this section, we bound the cop number of the graph $\Gamma = C(G, S)^\sigma$ as in \cref{Thm:CayleyTwistedCopNoBdd}. Note that the twisted Cayley graph $C(G, S)^\sigma$ is undirected if and only if $S$ contains $\sigma  (s^\mo g^\mo ) \sigma^\mo (g) $ for any $s\in S, g\in G$. 
Since $C(G, S)^\sigma$ is assumed to be undirected in \cref{Thm:CayleyTwistedCopNoBdd}, it follows that $S = \sigma(S^\mo)$ by taking $g = e$, and hence $\sigma(S) = S^\mo$. 
In the proof of \cref{Thm:CayleyTwistedCopNoBdd}, the assumption that $S$ is symmetric is used crucially, for instance, to make sure that the proposed strategy for the cops consists of moves that are possible in $C(G, S)^\sigma$. Moreover, since $\sigma(S) = S^\mo = S$ and $C(G, S)^\sigma$ is connected, it follows that $S$ generates $G$. 
\subsection{The moves of the primary cops}

Let $n$ be a positive integer. Suppose at the beginning of the $n$-th round, the robber is at the vertex $y$ and moves to $\tilde y =\sigma( yt)$ for some $t\in S$, and a cop $C$ having label $c(s)$ is positioned at the vertex $z$. The \emph{connection} of the robber is defined as 
$$
\xi_{y, \widetilde y} 
: = 
\begin{cases}
t & \text{ if $n$ is odd}, \\
\sigma(t) & \text{ if  $n$ is even}. 
\end{cases}
$$
The \emph{gap} between the cop $C$ and the robber at the beginning of the $n$-th round is defined as 
$$
\Delta_{z, y}:=
\begin{cases}
z^\mo y & \text{ if $n$ is odd},\\
\sigma(z^\mo y) & \text{ if $n$ is even}.
\end{cases}
$$
We describe the move of the cop $C$ with label $c(s)$, and the method of relabelling $C$. 
If $n$ is odd, then 
$C$ moves to 
$$
\begin{cases}
\sigma(zt) & \text{ when Case (1) holds}, \\
\sigma(z g_1) & \text{ when Case (2a) holds}, \\
\sigma(zt) & \text{ when Case (2bi) holds}, \\
\sigma(zs) & \text{ when Case (2bii) holds}.
\end{cases}
$$
If $n$ is even, then $C$ moves to 
$$
\begin{cases}
\sigma(zt) & \text{ when Case (1) holds}, \\
\sigma(z) g_1& \text{ when Case (2a) holds}, \\
\sigma(zt)& \text{ when Case (2bi) holds}, \\
\sigma(z) s& \text{ when Case (2bii) holds}.
\end{cases}
$$
Note that the above moves are possible in $C(G, S)^\sigma$ since $S = \sigma(S)$, which follows from the assumption that $S$ is symmetric and $\sigma(S) = S^\mo$. The label of $C$ is changed to $c(\xi_{y, \widetilde y}^\mo s \xi_{y, \widetilde y})$. Since $S$ is closed under conjugation, the element $\xi_{y, \widetilde y}^\mo s \xi_{y, \widetilde y}$ lies in $S$ and hence, the relabelling is justified.

\subsubsection{Reduction of tail and power}
We prove \cref{Lemma:ReductionGen} for $\Gamma = C(G, S)^\sigma$, showing that under the above strategy, the tail and the power of the cops do not increase, and in certain situations, they strictly decrease.

\begin{proof}
[Proof of \cref{Lemma:ReductionGen} for $\Gamma = C(G, S)^\sigma$]
Assume that at the beginning of the $n$-th round, the robber is at the vertex $y$, and a cop $C$ is at the vertex $z$ with label $c(s)$, and the robber moves to $\sigma(yt)$. Suppose the cop $C$ moves to $\tilde z$, and its label changes to $c(\tilde s)$.

Suppose Case (1) holds, that is, the connection of the robber in the $n$-th round is not a power of the label that the cop $C$ has at the beginning of the $n$-th round. 
If $n$ is odd, then writing $z^\mo y = w_{z, y, k} s^i$, we note that the gap between the cop $C$ and the robber at the beginning of the $(n+1)$-st round satisfies 
\begin{align*}
\Delta_{\widetilde z, \widetilde y}
& = 
\sigma(\tilde z^\mo \tilde y)
\\
& = 
\sigma ( \sigma(zt)^\mo \sigma(y t)) \\
& =
t^\mo z^\mo y t\\
& = 
(t^\mo w_{z, y, k} t) (t^\mo s t)^i\\
& = 
(t^\mo w_{z, y, k} t) (\xi_{y, \widetilde y}^\mo s \xi_{y, \widetilde y})^i\\
& = 
(t^\mo w_{z, y, k} t) \tilde s^i.
\end{align*}
If $n$ is even, then writing $\sigma(z^\mo y) = w_{z, y, k} s^i$, we note that the gap between the cop $C$ and the robber at the beginning of the $(n+1)$-st round satisfies 
\begin{align*}
\Delta_{\widetilde z, \widetilde y}
& = 
\tilde z^\mo \tilde y
\\
& = 
\sigma(zt)^\mo \sigma(y t) \\
& =
(\sigma(t)^\mo
 w_{z, y, k} \sigma(t))
 (\sigma(t)^\mo s \sigma(t))^i \\
& =
(\sigma(t)^\mo
 w_{z, y, k} \sigma(t))
( \xi_{y, \widetilde y}^\mo s \xi_{y, \widetilde y})^i \\
& =
(\sigma(t)^\mo
 w_{z, y, k} \sigma(t))
\tilde s^i 
.
\end{align*}
This proves the first statement of the lemma in Case (1). 

Suppose Case (2) holds, that is, the connection of the robber in the $n$-th round is a power of the label that the cop $C$ has at the beginning of the $n$-th round, that is, $\xi_{y, \widetilde y} = s^j$ for some integer $j$. 
Note that the label of $C$ remains unchanged in this round. 
First, let us consider Case (2a). So, the tail of the cop $C$ at the beginning of the $n$-th round is $k$ for some integer $k\geq 1$. 
If $n$ is odd, then writing $z^\mo y = g_1 g_2 \cdots g_k s^i$, we note that the gap between the cop $C$ and the robber at the beginning of the $(n+1)$-st round satisfies \begin{align*}
\Delta_{\widetilde z, \widetilde y}
& = 
\sigma(\tilde z^\mo \tilde y)
\\
& = \sigma(\sigma(z g_1)^\mo  \sigma(yt))\\
& = g_2 \cdots g_k s^i t\\
& = g_2 \cdots g_k s^i \xi_{y, \widetilde y}\\
& = g_2 \cdots g_k s^{i + j}
.
\end{align*}
If $n$ is even, then writing $\sigma(z^\mo y)= g_1 g_2 \cdots g_k s^i$, we note that the gap between the cop $C$ and the robber at the beginning of the $(n+1)$-st round satisfies 
\begin{align*}
\Delta_{\widetilde z, \widetilde y}
& = 
\tilde z^\mo \tilde y
\\
& = (\sigma(z) g_1)^\mo  \sigma(yt)\\
& = g_2 \cdots g_k s^i  \sigma(t)\\
& = g_2 \cdots g_k s^i \xi_{y, \widetilde y}\\
& = g_2 \cdots g_k s^{i + j}
.
\end{align*}
This proves the first statement in Case (2a), and also proves the second statement. 
Next, let us consider Case (2b). So, the tail of the cop $C$ at the beginning of the $n$-th round is zero. 
If $n$ is odd, then writing $z^\mo y = s^i$, we note that
\begin{align*}
\Delta_{\widetilde z, \widetilde y}
& = 
\sigma(\tilde z^\mo \tilde y)
\\
& = 
\sigma ( \sigma(zt)^\mo \sigma(y t)) \\
& =
t^\mo z^\mo y t\\
& =
s^i,
\end{align*}
and also note that 
\begin{align*}
\Delta_{\widetilde z, \widetilde y}
& = 
\sigma(\tilde z^\mo \tilde y)
\\
& = 
\sigma ( \sigma(zs)^\mo \sigma(y t)) \\
& =
s^\mo z^\mo y t\\
& = s^{i-1}t\\
& = s^{i-1} \xi_{y, \widetilde y},
\end{align*}
which is equal to $s^i$ (resp. $s^{i-2}$) when the connection $\xi_{y, \widetilde y}$ is equal to $s$ (resp. $s^\mo$). 
If $n$ is even, then writing $\sigma(z^\mo y) = s^i$, we note that 
\begin{align*}
\Delta_{\widetilde z, \widetilde y}
& = 
\tilde z^\mo \tilde y
\\
& = 
 \sigma(zt)^\mo \sigma(y t)\\
& =
(\sigma(t)^\mo s \sigma(t))^i
,
\end{align*}
and also note that 
\begin{align*}
\Delta_{\widetilde z, \widetilde y}
& = 
\tilde z^\mo \tilde y
\\
& = 
( \sigma(z)s)^\mo \sigma(y t) \\
& =
 s^{i-1}  \sigma(t)\\
& = s^{i-1} \xi_{y, \widetilde y}
,
\end{align*}
which is equal to $s^i$ (resp. $s^{i-2}$) when the connection $\xi_{y, \widetilde y}$ is equal to $s$ (resp. $s^\mo$). 
This proves the first statement in Case (2b), and also proves the third statement. This completes the proof of the lemma. 
\end{proof}

\subsection{The moves of the auxiliary cops} 
Let us introduce another set of $|S|$ \emph{auxiliary} cops and describe their moves. 
For each primary cop $C$ already employed, we introduce a cop $C_\aux$. During any round, $C_\aux$ occupies the vertex $\sigma(v)$ if $C$ is at the vertex $v$. Note that if $C$ moves from $u$ to $v$ during a round, then $v = \sigma(us)$ for some $s\in S$, and hence $\sigma(v) = \sigma( \sigma (u) t)$ with $t = \sigma(s)$ lying in $\sigma(S) = S^\mo = S$, and consequently, the cop $C_\aux$ can move from $\sigma(u)$ to $\sigma(v)$ in the graph $C(G, S)^\sigma$.

\subsection{Capturing the robber}
Using \cref{Lemma:ReductionGen}, we show that the cop number of $C(G, S)^\sigma$ is at most $2|S|$.

\begin{proof}
[Proof of \cref{Thm:CayleyTwistedCopNoBdd}]
\cref{Thm:CayleyTwistedCopNoBdd} can be established using \cref{Lemma:ReductionGen} for $\Gamma = C(G, S)^\sigma$ in a manner similar to the proof of 
\cref{Thm:CayleySumCopNoBdd}, which uses \cref{Lemma:ReductionGen} for $\Gamma = C_\Sigma(G, S)$. 

Indeed, writing the position $x$ of the robber at the end of the zeroth round as a product of $m$ elements of $S$, the pigeonhole principle implies that within at most $(m+ |G|)|S|$ rounds, there are precisely $m + |G|$ rounds (henceforth called the \emph{special rounds}), such that for any of these $m+|G|$ rounds, the connection of the robber is equal to the label that the cop $C$ has at the beginning of that round. Let $D$ denote the cop such that the labels of $C, D$ at the beginning of the first round (and hence, at the beginning of any round) are inverses of each other. 
Since the cops occupy the vertex at the identity element of $G$ and $x$ can be expressed as a product of $m$ elements of $G$, it follows that the tail of every cop at the beginning of the first round is at most $m$. 
By \cref{Lemma:ReductionGen}, the tails of the cops $C, D$ are equal to zero during the $(m+1)$-st special round and the subsequent rounds. 

If the power of $D$ is equal to $2\ell$ for some integer $0 \leq \ell \leq |G|/2$ at the beginning of the $(m+1)$-st special round, then by \cref{Lemma:ReductionGen}, its power is equal to zero at the end of the $(m+ \ell)$-th special round, i.e., the cop $D$ captures the robber at the end of the $(m+\ell)$-th special round. 
Suppose the power of  $D$ is equal  to $2\ell+1$ for some integer $\ell$  satisfying $1 \leq 2\ell+1 \leq |G|$ at the beginning of the $(m+1)$-st special round. By \cref{Lemma:ReductionGen}, the power of $D$ is equal to $1$ at the beginning of the $(m+\ell+1)$-st special round. Denote the position of the robber (resp. the cop $D$) at the beginning of the $(m+\ell+1)$-st special round by $y$ (resp. $z$). 
Let $s$ denote the label of $C$ at the beginning of this round. Suppose the robber moves to $\sigma(yt)$ during this round. 
If the $(m+\ell+1)$-st special round is an even round, then $\sigma(z^\mo y) = s^\mo$ and the connection of the robber is equal to $\sigma(t)$, which is equal to $s$, and consequently, $\sigma(yt) = \sigma(y) s = \sigma(z)$. 
If the $(m+\ell+1)$-st special round is an odd round, then $z^\mo y = s^\mo$ and the connection of the robber is equal to $t$, which is equal to $s$, and consequently, $\sigma(yt) = \sigma(ys) = \sigma(z)$. 
Hence, in this special round, the robber is captured by the auxiliary cop $D_\aux$. Consequently, the robber will be captured within at most $(m + |G|/2)|S|$ rounds. 
\end{proof}

\section{Cop numbers of twisted Cayley sum graphs}

Let $\sigma$ and $S$ be as in \cref{Thm:CayleyTwistedSumCopNoBdd}. 
In this section, we bound the cop number of the graph $\Gamma = C_\Sigma(G, S)^\sigma$ as in \cref{Thm:CayleyTwistedSumCopNoBdd}. 
Note that the twisted Cayley sum graph $C_\Sigma(G, S)^\sigma$ is undirected if and only if $S$ contains $\sigma^2  (g) \sigma(s) g^\mo$ for any $s\in S, g\in G$ \cite[Lemma 4.1]{CheegerTwisted}. 
Since $C_\Sigma(G, S)^\sigma$ is assumed to be undirected in \cref{Thm:CayleyTwistedSumCopNoBdd}, it follows that $S = \sigma(S)$ by taking $g = e$. 
Since $\sigma$ is an automorphism of order two and $\sigma(S) = S$, 
it follows that $S$ contains all conjugates of 
$\sigma(S)$, and as a consequence, $S$ is closed under conjugation. 
In the proof of \cref{Thm:CayleyTwistedSumCopNoBdd}, the assumption that $S$ is symmetric is used crucially, for instance, to make sure that the proposed strategy for the cops consists of moves that are possible in $C_\Sigma(G, S)^\sigma$. Moreover, since $\sigma(S) = S = S^\mo$ and $C_\Sigma(G, S)^\sigma$ is connected, it follows that $S$ generates $G$.

\subsection{The moves of the primary cops}

Let $n$ be a positive integer. Suppose at the beginning of the $n$-th round, the robber is at the vertex $y$ and moves to $\tilde y =\sigma( y^\mo t)$ for some $t\in S$, and a cop $C$ having label $c(s)$ is positioned at the vertex $z$. The \emph{connection} of the robber is defined as 
$$
\xi_{y, \widetilde y} 
: = 
\begin{cases}
y^\mo t^\mo y & \text{ if $n$ is odd} , \\
\sigma(t) & \text{ if $n$ is even} .
\end{cases}
$$
The \emph{gap} between the cop $C$ and the robber at the beginning of the $n$-th round is defined as 
$$
\Delta_{z, y}:=
\begin{cases}
z^\mo y & \text{ if $n$ is odd},\\
\sigma(z y^\mo) & \text{ if $n$ is even}.
\end{cases}
$$
We describe the move of the cop $C$ with label $c(s)$, and the method of relabelling $C$. 
If $n$ is odd, then 
$C$ moves to 
$$
\begin{cases}
\sigma(\xi_{y, \widetilde y}^\mo z^\mo) & \text{ when Case (1) holds}, \\
\sigma(g_1^\mo z^\mo) & \text{ when Case (2a) holds}, \\
\sigma(\xi_{y, \widetilde y}^\mo z^\mo) & \text{ when Case (2bi) holds}, \\
\sigma(s^\mo z^\mo) & \text{ when Case (2bii) holds}.
\end{cases}
$$
Note that the above moves are possible in $C_\Sigma(G, S)^\sigma$ since for any $x\in S$, there is an edge from $z$ to any vertex of the form $\sigma(x^\mo z^\mo)$ because $\sigma(x^\mo z^\mo )= \sigma(z^\mo (z x^\mo z^\mo))$ and $S$ is closed under conjugation and $S$ is symmetric. 
If $n$ is even, then $C$ moves to 
$$
\begin{cases}
\sigma(z^\mo t)& \text{ when Case (1) holds}, \\
\sigma(z^\mo) g_1& \text{ when Case (2a) holds}, \\
\sigma(z^\mo t)& \text{ when Case (2bi) holds}, \\
\sigma(z^\mo) s& \text{ when Case (2bii) holds}.
\end{cases}
$$
The label of $C$ is changed to $c(\xi_{y, \widetilde y}^\mo s \xi_{y, \widetilde y})$. Since $S$ is closed under conjugation, the element $\xi_{y, \widetilde y}^\mo s \xi_{y, \widetilde y}$ lies in $S$ and hence, the relabelling is justified.

\subsubsection{Reduction of tail and power}

We prove \cref{Lemma:ReductionGen} for $\Gamma = C_\Sigma(G, S)^\sigma$, showing that under the above strategy, the tail and the power of the cops do not increase, and in certain situations, they strictly decrease.

\begin{proof}
[Proof of \cref{Lemma:ReductionGen} for $\Gamma = C_\Sigma(G, S)^\sigma$]
Assume that at the beginning of the $n$-th round, the robber is at the vertex $y$, and a cop $C$ is at the vertex $z$ with label $c(s)$, and the robber moves to $\sigma(y^\mo t)$. 

Suppose Case (1) holds, that is, the connection of the robber in the $n$-th round is not a power of the label that the cop $C$ has at the beginning of the $n$-th round. 
If $n$ is odd, then writing $z^\mo y = w_{z, y, k} s^i$, we note that the gap between the cop $C$ and the robber at the beginning of the $(n+1)$-st round satisfies 
\begin{align*}
\Delta_{\widetilde z, \widetilde y}
& = 
\sigma(\tilde z \tilde y^\mo)
\\
& =
\sigma
\left(
\sigma(\xi_{y, \widetilde y}^\mo z^\mo)
(\sigma(y^\mo t))^\mo
\right) \\
& = \xi_{y, \widetilde y}^\mo z^\mo y (y^\mo t^\mo y) \\
& = \xi_{y, \widetilde y}^\mo w_{z, y, k} s^i \xi_{y, \widetilde y}\\
& = (\xi_{y, \widetilde y}^\mo w_{z, y, k} \xi_{y, \widetilde y} )
(\xi_{y, \widetilde y}^\mo s \xi_{y, \widetilde y} )^i\\
& = (\xi_{y, \widetilde y}^\mo 
w_{z, y, k}\xi_{y, \widetilde y})
\tilde s^i .
\end{align*}
If $n$ is even, then writing $\sigma(zy ^\mo) = w_{z, y, k} s^i$, we note that the gap between the cop $C$ and the robber at the beginning of the $(n+1)$-st round satisfies 
\begin{align*}
\Delta_{\widetilde z, \widetilde y}
& = 
\tilde z^\mo \tilde y
\\
& = (\sigma(z^\mo t) )^\mo (\sigma(y^\mo t) )\\
& = \sigma(t)^\mo \sigma(zy^\mo )\sigma(t )\\
& =  \sigma(t)^\mo  w_{z, y, k} s^i \sigma(t )\\
& =  (\sigma(t)^\mo  w_{z, y, k} \sigma(t) )(\sigma(t)^\mo  s \sigma(t) )^i \\
& =  (\sigma(t)^\mo  w_{z, y, k} \sigma(t) )(\xi_{y, \widetilde y}^\mo  s \xi_{y, \widetilde y})^i \\
& =  (\sigma(t)^\mo  w_{z, y, k} \sigma(t) )
\tilde s^i 
.
\end{align*}
This proves the first statement of the lemma in Case (1). 

Suppose Case (2) holds, that is, the connection of the robber in the $n$-th round is a power of the label that the cop $C$ has at the beginning of the $n$-th round, that is, $\xi_{y, \widetilde y} = s^j$ for some integer $j$. 
First, let us consider Case (2a). So, the tail of the cop $C$ at the beginning of the $n$-th round is $k$ for some integer $k\geq 1$. 
If $n$ is odd, then writing $z^\mo y = g_1 g_2 \cdots g_k s^i$, we note that the gap between the cop $C$ and the robber at the beginning of the $(n+1)$-st round satisfies \begin{align*}
\Delta_{\widetilde z, \widetilde y}
& = 
\sigma(\tilde z \tilde y^\mo)
\\
& =
\sigma(\sigma(g_1^\mo z^\mo)  (\sigma(y^\mo t))^\mo )\\
& = g_2 \cdots g_k s^i y^\mo t^\mo y \\
& = g_2 \cdots g_k s^i \xi_{y, \widetilde y} \\
& = g_2 \cdots g_k s^{i + j}
.
\end{align*}
If $n$ is even, then writing $\sigma(zy^\mo)= g_1 g_2 \cdots g_k s^i$, we note that the gap between the cop $C$ and the robber at the beginning of the $(n+1)$-st round satisfies 
\begin{align*}
\Delta_{\widetilde z, \widetilde y}
& = 
\tilde z^\mo \tilde y
\\
& =
(\sigma(z)^\mo g_1)^\mo \sigma(y^\mo t)\\
& = g_1^\mo \sigma(z y^\mo) \sigma(t)\\
& = g_2 \cdots g_k s^i  \sigma(t)\\
& = g_2 \cdots g_k s^i \xi_{y, \widetilde y} \\
& = g_2 \cdots g_k s^{i + j}
.
\end{align*}
This proves the first statement in Case (2a), and also proves the second statement. 
Next, let us consider Case (2b). So, the tail of the cop $C$ at the beginning of the $n$-th round is zero. 
If $n$ is odd, then writing $z^\mo y = s^i$, we note that
\begin{align*}
\Delta_{\widetilde z, \widetilde y}
& = 
\sigma(\tilde z \tilde y^\mo)
\\
& =
\sigma(\sigma(\xi_{y, \widetilde y}^\mo z^\mo )\sigma(  (y^\mo t))^\mo ) \\
& = 
\sigma(\sigma((y^\mo t^\mo y)^\mo z^\mo )\sigma(  (y^\mo t))^\mo ) \\
& = (y^\mo t^\mo y)^\mo s^{i} y^\mo  t^\mo y\\
& = s^i,
\end{align*}
and also note that 
\begin{align*}
\Delta_{\widetilde z, \widetilde y}
& = 
\sigma(\tilde z \tilde y^\mo)
\\
& =
\sigma(\sigma(s^\mo z^\mo)  (\sigma(y^\mo t))^\mo )\\
& = \sigma(\sigma(s^\mo ) \sigma(s^i y^\mo )\sigma( t^\mo y))\\
& = s^{i-1} y^\mo  t^\mo y \\
& = s^{i-1} \xi_{y, \widetilde y}
,
\end{align*}
which is equal to $s^i$ (resp. $s^{i-2}$) when the connection $\xi_{y, \widetilde y}$ is equal to $s$ (resp. $s^\mo$). 
If $n$ is even, then writing $\sigma(z y^\mo) = s^i$, we note that 
\begin{align*}
\Delta_{\widetilde z, \widetilde y}
& = 
\tilde z^\mo \tilde y
\\
& =
(\sigma(z^\mo t))^\mo \sigma(y^\mo t)\\
& = \sigma(t)^\mo \sigma(zy^\mo ) \sigma(t)\\
& = \sigma(t)^\mo s^{i} \sigma(t)\\
& = s^i
,
\end{align*}
and also note that 
\begin{align*}
\Delta_{\widetilde z, \widetilde y}
& = 
\tilde z^\mo \tilde y
\\
& = (\sigma(z)^\mo s)^\mo\sigma( y^\mo t)\\
& = s^\mo \sigma(zy^\mo )\sigma( t)\\
& = s^{i-1} \sigma(t) \\
& = s^{i-1} \xi_{y, \widetilde y}
,
\end{align*}
which is equal to $s^i$ (resp. $s^{i-2}$) when the connection $\xi_{y, \widetilde y}$ is equal to $s$ (resp. $s^\mo$). 
This proves the first statement in Case (2b), and also proves the third statement. This completes the proof of the lemma. 
\end{proof}

\subsection{The moves of the auxiliary cops} 
Let us introduce another set of $|S|$ \emph{auxiliary} cops and describe their moves. 
For each primary cop $C$ already employed, we introduce a cop $C_\aux$. During any round, $C_\aux$ occupies the vertex $\sigma(v^\mo)$ if $C$ is at the vertex $v$. 
Note that if $C$ moves from $u$ to $v$ during a round, then $u\sigma^\mo (v) = u \sigma(v)$ lies in the symmetric set $S$, and hence there is an edge joining the vertex $\sigma(v^\mo)$ with $\sigma(u^\mo)$ in $C_\Sigma(G, S)^\sigma$, and consequently, the cop $C_\aux$ can move from $\sigma(u^\mo)$ to $\sigma(v^\mo)$ in the undirected graph $C_\Sigma(G, S)^\sigma$.

\subsection{Capturing the robber}
Using \cref{Lemma:ReductionGen}, we show that the cop number of $C_\Sigma(G, S)^\sigma$ is at most $2|S|$.

\begin{proof}
[Proof of \cref{Thm:CayleyTwistedSumCopNoBdd}]
\cref{Thm:CayleyTwistedSumCopNoBdd} can be established using \cref{Lemma:ReductionGen} for $\Gamma = C_\Sigma(G, S)^\sigma$ in a manner similar to the proof of 
\cref{Thm:CayleySumCopNoBdd}, which uses \cref{Lemma:ReductionGen} for $\Gamma = C_\Sigma(G, S)$. 

Indeed, writing the position $x$ of the robber at the end of the zeroth round as a product of $m$ elements of $S$, the pigeonhole principle implies that within at most $(m+ |G|)|S|$ rounds, there are precisely $m + |G|$ rounds (henceforth called the \emph{special rounds}), such that for any of these $m+|G|$ rounds, the connection of the robber is equal to the label that the cop $C$ has at the beginning of that round. Let $D$ denote the cop such that the labels of $C, D$ at the beginning of the first round (and hence, at the beginning of any round) are inverses of each other. 
Since the cops occupy the vertex at the identity element of $G$ and $x$ can be expressed as a product of $m$ elements of $G$, it follows that the tail of every cop at the beginning of the first round is at most $m$. By \cref{Lemma:ReductionGen}, the tails of the cops $C, D$ are equal to zero during the $(m+1)$-st special round and the subsequent rounds. 

If the power of $D$ is equal to $2\ell$ for some integer $0 \leq \ell \leq |G|/2$ at the beginning of the $(m+1)$-st special round, then by \cref{Lemma:ReductionGen}, its power is equal to zero at the end of the $(m+ \ell)$-th special round, i.e., the cop $D$ captures the robber at the end of the $(m+\ell)$-th special round. 
Suppose the power of  $D$ is equal  to $2\ell+1$ for some integer $\ell$  satisfying $1 \leq 2\ell+1 \leq |G|$ at the beginning of the $(m+1)$-st special round. By \cref{Lemma:ReductionGen}, the power of $D$ is equal to $1$ at the beginning of the $(m+\ell+1)$-st special round. Denote the position of the robber (resp. the cop $D$) at the beginning of the $(m+\ell+1)$-st special round by $y$ (resp. $z$). 
Let $s$ denote the label of $C$ at the beginning of this round. Suppose the robber moves to $\sigma(y^\mo t)$ during this round. 
If the $(m+\ell+1)$-st special round is an even round, then $\sigma(zy^\mo) = s^\mo$ and the connection of the robber is equal to $\sigma(t)$, which is equal to $s$, and consequently, $\sigma(y^\mo t) = \sigma(y^\mo ) s = \sigma(z^\mo)$. 
If the $(m+\ell+1)$-st special round is an odd round, then $z^\mo y = s^\mo$ and the connection of the robber is equal to $y^\mo t^\mo y$, which is equal to $s$, and consequently, $\sigma(y^\mo t) = \sigma(s^\mo y^\mo)  = \sigma(z^\mo)$. 
Hence, in this special round, the robber is captured by the auxiliary cop $D_\aux$. Consequently, the robber will be captured within at most $(m + |G|/2)|S|$ rounds. 
\end{proof}

\section{The weak Meyniel conjecture}

The weak Meyniel conjecture states that for a fixed constant $c>0$, the cop number of a graph on $n$ vertices is $O(n^{1-c})$ 
\cite[p. 227]{BairdBonatoMeynielConjSurvey}, 
\cite[p. 33]{BradshawHosseiniMoharStachoOnTheCopNumberGraphWithHighGirth}. 
Bradshaw--Hosseini--Mohar--Stacho showed that the weak Meyniel conjecture holds for expander families of bounded degree 
\cite[Corollary 4.7]{BradshawHosseiniMoharStachoOnTheCopNumberGraphWithHighGirth}. In this section, under certain hypotheses, we prove a stronger version of the weak Meyniel conjecture for 
Cayley graphs, Cayley sum graphs, and twisted Cayley and Cayley sum graphs. 
We use a result of Bollob\'{a}s--Janson--Riordan to construct a ``small'' dominating set for these graphs when the degree is ``large'', and we use Theorems 
\ref{Thm:CayleySumCopNoBdd}, \ref{Thm:CayleyTwistedCopNoBdd}, \ref{Thm:CayleyTwistedSumCopNoBdd} when the degree is ``small''.

\begin{theorem}
\label{Thm:WMC}
\quad 
\begin{enumerate}
\item 
The cop number of any undirected, connected Cayley graph $C(G, S)$ of a group $G$ of order $n$ is at most $2\sqrt n \log n$, where $S$ is a subset of $G$, closed under conjugation by the elements of $G$.  
 
\item 

The cop number of any undirected, connected Cayley sum graph $C_\Sigma(G, S)$ of a group $G$ of order $n$ is at most $2\sqrt n \log n$, where $S$ is a symmetric subset of $G$. 
 
\item 
 
The cop number of any undirected, connected twisted Cayley graph $C(G, S)^\sigma$ of a group $G$ of order $n$ is at most $2\sqrt n \log n$, where $S$ is a symmetric subset of $G$, closed under conjugation by the elements of $G$, and $\sigma$ denotes an automorphism of $G$ of order two.

\item 

The cop number of any undirected, connected twisted Cayley sum graph $C_\Sigma(G, S)^\sigma$ of a group $G$ of order $n$ is at most $2\sqrt n \log n$, where $S$ is a symmetric subset of $G$ and $\sigma$ denotes an automorphism of $G$ of order two. 

\end{enumerate} 
Consequently, for any real number $c < \frac 12$, the above cop numbers are  $O( n^{1-c})$. 
\end{theorem}

\begin{proof}
If $S$ is a subset of $G$ satisfying $|S| \geq \sqrt n \log n$, then by \cite[Corollary 3.2]{BollobasJansonRiordanOnCoveringByTranslatesOfASet}, there exists a subset $T$ of $G$ of order $\leq \sqrt n (1 + \frac 1 {\log n}) \leq 2\sqrt n$ such that $TS = G$, i.e., $G$ can be covered by at most $2\sqrt n$ left translates of $S$. 

Note that $TS =G$ implies that $T$ is a dominating set of $C(G, S)$. By employing one cop at each of the vertices of $C(G, S)$ corresponding to the elements of $T$ when $|S| \geq \sqrt n \log n$, it follows that the robber will be captured in $C(G, S)$. 
Further, if $S$  is closed under conjugation by the elements of $G$,  $C(G, S)$ is connected and undirected, and $|S| < \sqrt n \log n$, then by \cite[Theorem 1.5]{FranklCopsRobbersGraphsLargeGirthCayleyGraphs}, the cop number of $C(G, S)$ is at most $\sqrt n \log n$. 

Note that $TS =G$ shows that $T^\mo$ is a dominating set of $C_\Sigma(G, S)$. 
By employing one cop at each of the vertices of $C_\Sigma(G, S)$ corresponding to the elements of $T^\mo$ when $|S| \geq \sqrt n \log n$, it follows that the robber will be captured in $C_\Sigma(G, S)$. 
Further, if $S$  is symmetric, $C_\Sigma(G, S)$ is connected and undirected, and $|S| < \sqrt n \log n$, then by applying \cref{Thm:CayleySumCopNoBdd}, it follows that the cop number of $C_\Sigma(G, S)$  is at most $2\sqrt n \log n$.

Note that $TS = G$ together with $\sigma(S) = S$ implies that $\sigma(T)$ is a dominating set of $C(G, S)^\sigma$. By employing one cop at each of the vertices of $C(G, S)^\sigma$ corresponding to the elements of $\sigma(T)$ when $|S| \geq \sqrt n \log n$, it follows that the robber will be captured in $C(G, S)^\sigma$. 
Further, if $S$  is  symmetric and it is closed under conjugation by the elements of $G$, and $C(G, S)^\sigma$ is connected and undirected, and $|S| < \sqrt n \log n$, then by applying \cref{Thm:CayleyTwistedCopNoBdd}, it follows that the cop number of $C(G, S)^\sigma$  is at most $2\sqrt n \log n$.

Note that $TS =G$ together with $\sigma(S) = S$ shows that $\sigma(T^\mo)$ is a dominating set of $C_\Sigma(G, S)^\sigma$. 
By employing one cop at each of the vertices of $C_\Sigma(G, S)^\sigma$ corresponding to the elements of $\sigma(T^\mo)$ when $|S| \geq \sqrt n \log n$, it follows that the robber will be captured in $C_\Sigma(G, S)^\sigma$. 
Further, if $S$  is symmetric, $C_\Sigma(G, S)^\sigma$ is connected and undirected, and $|S| < \sqrt n \log n$, then by applying \cref{Thm:CayleyTwistedSumCopNoBdd}, it follows that the cop number of $C_\Sigma(G, S)^\sigma$  is at most $2\sqrt n \log n$.
\end{proof}

\section{Acknowledgements}
The authors would like to thank the referees for providing constructive suggestions, which helped improve the manuscript. 
A.B. first learnt about the game of cops and robbers from Arnab Chakraborty at the Indian Statistical Institute, Kolkata in 2008, and wishes to thank him and also Kishalaya Saha for a number of helpful discussions on combinatorial games during that period. This work started during the first author's stay at the University of Copenhagen, where he was supported by the ERC grant 716424 - CASe of Karim Adiprasito. The second author acknowledges the INSPIRE Faculty Award IFA18-MA123 from the Department of Science and Technology, Government of India. 

\section*{Data availability statement}
Data sharing not applicable to this article as no datasets were generated or analysed during the current study. 

\section*{Conflict of interest} The authors declare that they have no conflict of interest.

\bibliographystyle{amsplain}

\begin{thebibliography}{10}

\bibitem{AignerFrommeGameOfCopsRobbers}
M.~Aigner and M.~Fromme, \emph{A game of cops and robbers}, Discrete Appl.
  Math. \textbf{8} (1984), no.~1, 1--11. \MR{739593}

\bibitem{AmooshahiTaeri}
Marzieh Amooshahi and Bijan Taeri, \emph{On {C}ayley sum graphs of non-abelian
  groups}, Graphs Combin. \textbf{32} (2016), no.~1, 17--29. \MR{3436946}

\bibitem{AndreaeNoteOnPursuitGamePlayedGraphs}
Thomas Andreae, \emph{Note on a pursuit game played on graphs}, Discrete Appl.
  Math. \textbf{9} (1984), no.~2, 111--115. \MR{761595}

\bibitem{AndreaePursuitGameGraphsMinorExcluded}
\bysame, \emph{On a pursuit game played on graphs for which a minor is
  excluded}, J. Combin. Theory Ser. B \textbf{41} (1986), no.~1, 37--47.
  \MR{854602}

\bibitem{BairdBonatoMeynielConjSurvey}
William Baird and Anthony Bonato, \emph{Meyniel's conjecture on the cop number:
  a survey}, J. Comb. \textbf{3} (2012), no.~2, 225--238. \MR{2980752}

\bibitem{CheegerCayleySum}
Arindam Biswas and Jyoti~Prakash Saha, \emph{A {C}heeger type inequality in
  finite {C}ayley sum graphs}, Algebr. Comb. \textbf{4} (2021), no.~3,
  517--531. \MR{4275826}

\bibitem{Expansion}
\bysame, \emph{Expansion in {C}ayley graphs, {C}ayley sum graphs and their
  twists}, Available at \url{https://arxiv.org/abs/2103.05935}, 2021.

\bibitem{CheegerTwisted}
\bysame, \emph{Spectra of twists of {C}ayley and {C}ayley sum graphs}, Adv. in
  Appl. Math. \textbf{132} (2022), Paper No. 102272, 34. \MR{4327334}

\bibitem{VertexTra}
\bysame, \emph{A spectral bound for vertex-transitive graphs and their spanning
  subgraphs}, Algebr. Comb. \textbf{6} (2023), no.~3, 689--706.

\bibitem{BollobasJansonRiordanOnCoveringByTranslatesOfASet}
B\'{e}la Bollob\'{a}s, Svante Janson, and Oliver Riordan, \emph{On covering by
  translates of a set}, Random Structures Algorithms \textbf{38} (2011),
  no.~1-2, 33--67. \MR{2768883}

\bibitem{BollobasKunLeaderCopsRobbersInRandomGraph}
B\'{e}la Bollob\'{a}s, G\'{a}bor Kun, and Imre Leader, \emph{Cops and robbers
  in a random graph}, J. Combin. Theory Ser. B \textbf{103} (2013), no.~2,
  226--236. \MR{3018067}

\bibitem{BonatoNowakowskiSTMLCopsRobbers}
Anthony Bonato and Richard~J. Nowakowski, \emph{The game of cops and robbers on
  graphs}, Student Mathematical Library, vol.~61, American Mathematical
  Society, Providence, RI, 2011. \MR{2830217}

\bibitem{BradshawProofOfMeynielConjAbelianCayleyGraphs}
Peter Bradshaw, \emph{A proof of the {M}eyniel conjecture for {A}belian
  {C}ayley graphs}, Discrete Math. \textbf{343} (2020), no.~1, 111546, 5.
  \MR{4039389}

\bibitem{BradshawHosseiniMoharStachoOnTheCopNumberGraphWithHighGirth}
Peter Bradshaw, Seyyed~Aliasghar Hosseini, Bojan Mohar, and Ladislav Stacho,
  \emph{On the cop number of graphs of high girth}, J. Graph Theory
  \textbf{102} (2023), no.~1, 15--34. \MR{4520060}

\bibitem{BradshawHosseiniTurcotteCopsRobbersDirUnditAbelianCayley}
Peter Bradshaw, Seyyed~Aliasghar Hosseini, and J\'{e}r\'{e}mie Turcotte,
  \emph{Cops and robbers on directed and undirected abelian {C}ayley graphs},
  European J. Combin. \textbf{97} (2021), Paper No. 103383, 19. \MR{4292350}

\bibitem{ChiniforooshanBetterBddForCopNumberGeneralGraph}
Ehsan Chiniforooshan, \emph{A better bound for the cop number of general
  graphs}, J. Graph Theory \textbf{58} (2008), no.~1, 45--48. \MR{2404040}

\bibitem{Chung89JAMS}
F.~R.~K. Chung, \emph{Diameters and eigenvalues}, J. Amer. Math. Soc.
  \textbf{2} (1989), no.~2, 187--196. \MR{965008}

\bibitem{DasGahlawatSahooSenCopsRobberSomeFamiliesOriGra}
Sandip Das, Harmender Gahlawat, Uma~kant Sahoo, and Sagnik Sen, \emph{Cops and
  {R}obber on some families of oriented graphs}, Theoret. Comput. Sci.
  \textbf{888} (2021), 31--40. \MR{4316755}

\bibitem{DeVosGoddynMoharSamalCayleySumFullerene}
Matt DeVos, Luis Goddyn, Bojan Mohar, and Robert {\v{S}}\'{a}mal, \emph{Cayley
  sum graphs and eigenvalues of {$(3,6)$}-fullerenes}, J. Combin. Theory Ser. B
  \textbf{99} (2009), no.~2, 358--369. \MR{2482954}

\bibitem{DudeneyAmusement}
Henry~Ernest Dudeney, \emph{Amusements in {M}athematics}, 1917.

\bibitem{FranklPursuitGameOnCayleyGraph}
P.~Frankl, \emph{On a pursuit game on {C}ayley graphs}, Combinatorica
  \textbf{7} (1987), no.~1, 67--70. \MR{905152}

\bibitem{FranklCopsRobbersGraphsLargeGirthCayleyGraphs}
Peter Frankl, \emph{Cops and robbers in graphs with large girth and {C}ayley
  graphs}, Discrete Appl. Math. \textbf{17} (1987), no.~3, 301--305.
  \MR{890640}

\bibitem{FriezeKrivelevichLohVariationsOnCopsRobbers}
Alan Frieze, Michael Krivelevich, and Po-Shen Loh, \emph{Variations on cops and
  robbers}, J. Graph Theory \textbf{69} (2012), no.~4, 383--402. \MR{2979296}

\bibitem{GonzalezSebastianHosseiniKnoxMoharReedCopsRobberOriToroiGrid}
Sebasti\'{a}n Gonz\'{a}lez Hermosillo de~la Maza, Seyyed~Aliasghar Hosseini,
  Fiachra Knox, Bojan Mohar, and Bruce Reed, \emph{Cops and robbers on oriented
  toroidal grids}, Theoret. Comput. Sci. \textbf{857} (2021), 166--176.
  \MR{4204531}

\bibitem{GreenCountingSetsWithSmallSubsetsClique}
Ben Green, \emph{Counting sets with small sumset, and the clique number of
  random {C}ayley graphs}, Combinatorica \textbf{25} (2005), no.~3, 307--326.
  \MR{2141661}

\bibitem{GreenChromatic}
\bysame, \emph{On the chromatic number of random {C}ayley graphs}, Combin.
  Probab. Comput. \textbf{26} (2017), no.~2, 248--266. \MR{3603967}

\bibitem{GreenMorrisCountingSetsWithSmallSubsets}
Ben Green and Robert Morris, \emph{Counting sets with small sumset and
  applications}, Combinatorica \textbf{36} (2016), no.~2, 129--159.
  \MR{3516881}

\bibitem{GrynkiewiczLevSerraConnCaylSum}
David Grynkiewicz, Vsevolod~F. Lev, and Oriol Serra, \emph{Connectivity of
  addition {C}ayley graphs}, J. Combin. Theory Ser. B \textbf{99} (2009),
  no.~1, 202--217. \MR{2467826}

\bibitem{HamidounePursuitGameOnCayleyGraphs}
Yahya~Ould Hamidoune, \emph{On a pursuit game on {C}ayley digraphs}, European
  J. Combin. \textbf{8} (1987), no.~3, 289--295. \MR{919881}

\bibitem{HasiriShinkarMeynielExtremalFamiliesAbelianCayley}
Fatemeh Hasiri and Igor Shinkar, \emph{Meyniel extremal families of abelian
  {C}ayley graphs}, Graphs Combin. \textbf{38} (2022), no.~3, Paper No. 61, 13.
  \MR{4393986}

\bibitem{HosseiniMoharGonzalezMeynielConjGraphsBddDegre}
Seyyed~Aliasghar Hosseini, Bojan Mohar, and Sebastian Gonzalez Hermosillo de~la
  Maza, \emph{Meyniel's conjecture on graphs of bounded degree}, J. Graph
  Theory \textbf{97} (2021), no.~3, 401--407. \MR{4313187}

\bibitem{LehnerCopNumberOfToroidalGraphs}
Florian Lehner, \emph{On the cop number of toroidal graphs}, J. Combin. Theory
  Ser. B \textbf{151} (2021), 250--262. \MR{4285899}

\bibitem{LevSumDiffHamiltCycle}
Vsevolod~F. Lev, \emph{Sums and differences along {H}amiltonian cycles},
  Discrete Math. \textbf{310} (2010), no.~3, 575--584. \MR{2564813}

\bibitem{LuPengMeynielConjOfTheCopNumber}
Linyuan Lu and Xing Peng, \emph{On {M}eyniel's conjecture of the cop number},
  J. Graph Theory \textbf{71} (2012), no.~2, 192--205. \MR{2965383}

\bibitem{MaamounMeynielGameOfPolicemenRobber}
M.~Maamoun and H.~Meyniel, \emph{On a game of policemen and robber}, Discrete
  Appl. Math. \textbf{17} (1987), no.~3, 307--309. \MR{890641}

\bibitem{MarusicScapellatoZagagliaSalviGeneralizedCayleyGraph}
Dragan Maru\v{s}i\v{c}, Raffaele Scapellato, and Norma Zagaglia~Salvi,
  \emph{Generalized {C}ayley graphs}, Discrete Math. \textbf{102} (1992),
  no.~3, 279--285. \MR{1169147}

\bibitem{NowakowskiWinklerVertexToVertexPursuitInGraph}
Richard Nowakowski and Peter Winkler, \emph{Vertex-to-vertex pursuit in a
  graph}, Discrete Math. \textbf{43} (1983), no.~2-3, 235--239. \MR{685631}

\bibitem{PralatWhenDoesARandomGraphHaveConsCopNumber}
Pawe{\l } Pra{\l}at, \emph{When does a random graph have constant cop number?},
  Australas. J. Combin. \textbf{46} (2010), 285--296. \MR{2598712}

\bibitem{PralatWormaldMeynielConjRandomGraphs}
Pawe{\l } Pra{\l}at and Nicholas Wormald, \emph{Meyniel's conjecture holds for
  random graphs}, Random Structures Algorithms \textbf{48} (2016), no.~2,
  396--421. \MR{3449604}

\bibitem{PralatWormaldMeynielConjRandomDReguGraphs}
\bysame, \emph{Meyniel's conjecture holds for random {$d$}-regular graphs},
  Random Structures Algorithms \textbf{55} (2019), no.~3, 719--741.
  \MR{3997485}

\bibitem{QuilliotThesis}
Alain Quilliot, \emph{Jeux et pointes fixes sur les graphes}, 1978, Th\`{e}se
  de 3\`{e}me cycle, Universit\'{e} de Paris VI.

\bibitem{QuilliotDiscretePursuitGame}
\bysame, \emph{Discrete pursuit game}, Congr. Numer. \textbf{38} (1983),
  227--241. \MR{703252}

\bibitem{ScottSudakovBddForTheCopsRobbersProblem}
Alex Scott and Benny Sudakov, \emph{A bound for the cops and robbers problem},
  SIAM J. Discrete Math. \textbf{25} (2011), no.~3, 1438--1442. \MR{2837608}

\end{thebibliography}
\end{document}